\documentclass[10pt]{article}
\usepackage{bm,latexsym,amsmath,amsthm,amssymb}
\usepackage{graphicx,indentfirst}
\usepackage[numbers,sort&compress]{natbib}
\date{}
\begin{document}
\title{\bf Continuity properties of the data-to-solution map for the two-component higher order Camassa-Holm system}
\author{Feng Wang\thanks{Corresponding author.\newline
\mbox{}\qquad E-mail: wangfeng@xidian.edu.cn (F. Wang); fqli@dlut.edu.cn (F. Li)}\\
\small School of Mathematics and Statistics, Xidian University,
    Xi'an 710071, PR China\\[3pt]
Fengquan Li\\
\small School of Mathematical Sciences, Dalian University of Technology,
Dalian 116024, PR China}
\date{}
\maketitle \baselineskip 3pt
\begin{center}
\begin{minipage}{130mm}
{{\bf Abstract.} This work studies the Cauchy problem of a two-component higher order Camassa-Holm system, which is well-posed in Sobolev spaces $H^{s}(\mathbb{R})\times H^{s-2}(\mathbb{R})$, $s>\frac{7}{2}$ and its solution map is continuous. We show that the solution map is H\"{o}lder continuous in $H^{s}(\mathbb{R})\times H^{s-2}(\mathbb{R})$ equipped with the $H^{r}(\mathbb{R})\times H^{r-2}(\mathbb{R})$-topology for $1\leq r<s$, and the H\"{o}lder exponent is expressed in terms of $s$ and $r$.

}

\end{minipage}
\end{center}

\baselineskip=15pt

\section{Introduction}
\label{intro}

In this paper, we consider the Cauchy problem of the following two-component higher order Camassa-Holm system
$$
\begin{array}{l}
\left\{\begin{array}{l}
m_{t}=\alpha u_{x}-bu_{x}m-um_{x}-\kappa \rho\rho_{x}, \quad m=Au, \\[3pt]
\rho_{t}=-u\rho_{x}-(b-1)u_{x}\rho, \quad b\in \mathbb{R}\setminus \{1\},\\[3pt]
\alpha_{t}=0,
\end{array}
\right.
\end{array}
\eqno(1.1)
$$
where $Au=(1-\partial_{x}^{2})^{\sigma}u$ with $\sigma>1$, and $b, \kappa$ are real parameters. Eq.(1.1) was proposed by Escher and Lyons \cite{esly15}, in which they showed that the system corresponds to a metric induced geodesic flow on the infinite dimensional Lie group $\textmd{Diff}^{\infty}(\mathbb{S}^{1})\circledS C^{\infty}(\mathbb{S}^{1})\times \mathbb{R}$ and admits a global solution in $C^{\infty}([0, \infty); C^{\infty}(\mathbb{S}^{1})\oplus C^{\infty}(\mathbb{S}^{1}))$ with smooth initial data in $C^{\infty}(\mathbb{S}^{1})\oplus C^{\infty}(\mathbb{S}^{1}))$ when $b=2$, where $\textmd{Diff}^{\infty}(\mathbb{S}^{1})$ denotes the group of orientation preserving diffeomorphisms of the circle and $\circledS$ denotes an appropriate semi-direct product between the pair. Recently, He and Yin \cite{heyin14}, Chen and Zhou \cite{chzh17} established the local well-posedness of (1.1) in Besov spaces. Zhou \cite{zhou18}, Zhang and Li \cite{zhli17} investigated the local well-posedness, blow-up criteria and Gevrey regularity of the solutions to (1.1) with $\sigma=2$. When $\rho\equiv 0$, $\alpha=0$ and $b=2$, (1.1) reduces to a Camassa-Holm equation with fractional order inertia operator, whose geometrical interpretation and local well-posedness can be seen in \cite{esko14,eskol14,heyin14}, and if we further assume $2\leq\sigma\in \mathbb{Z}_{+}$, (1.1) becomes a higher order Camassa-Holm equation derived as the Euler-Poincar\'{e} differential equation on the Bott-Virasoro group with respect to the $H^{\sigma}$ metric \cite{mczh09}.

For $\sigma=1$, (1.1) reduces to the following nonlinear system \cite{ehkl16}
$$
\begin{array}{l}
\left\{\begin{array}{l}
m_{t}=\alpha u_{x}-bu_{x}m-um_{x}-\kappa \rho\rho_{x}, \quad m=u-u_{xx}, \\[3pt]
\rho_{t}=-u\rho_{x}-(b-1)u_{x}\rho, \quad b\in \mathbb{R}\setminus \{1\},
\end{array}
\right.
\end{array}
\eqno(1.2)
$$
which models the two-component shallow water waves with constant vorticity $\alpha$. In \cite{ehkl16}, Escher et al. showed the local well-posedness of (1.2) under a geometrical framework, and studied the blow-up scenarios and global strong solutions of (1.2) on the circle. In \cite{guhey15}, Guan et al. considered the Cauchy problem of (1.2) in the Besov space and showed that the solutions have exponential decay if the initial data has exponential decay. When $\alpha=0$, $b=2$ and $\kappa=\pm 1$, (1.2) becomes the two-component Camassa-Holm system, which admits Lax pair
and bi-Hamiltonian structure, and thus is completely integrable \cite{chlz06}. When $\rho\equiv 0$ and $\alpha=0$, (1.2) reduces to a family of equations parameterised by $b\neq 1$, the so-called $b$-family equation. In particular, when $b=2$ and $b=3$, the $b$-family equation respectively becomes the famous completely integrable Camassa-Holm equation \cite{ch93} and Degasperis-Procesi equation \cite{desp99}, which were introduced to model the unidirectional propagation of shallow water waves over a flat bottom. The Cauchy problem for these equations have been well-studied both on the real line and on the circle, including the well-posedness, blow-up behavior, global existence, traveling wave solutions and so on, e.g. \cite{brc07,ce98,ces98,cm99,cm00,elyi07,esyi08,eslyin06,guiyin10,guiyin11,
guil10,himkm10,hmkz07,lene05a,lene05b,lo00,xz00} and the references
therein.

The present paper is devoted to establishing the H\"{o}lder continuity of the data-to-solution map for system (1.1) with $\sigma=2$ in $H^{s}(\mathbb{R})\times H^{s-2}(\mathbb{R})$, $s>\frac{7}{2}$, which provides more information about the stability of the solution map than the one given by Corollary 3.1.2 in \cite{zhli17}. We mention that H\"{o}lder continuity for the $b$-equation was proved on the line by Chen, Liu and Zhang in \cite{chlz13}, and for other equations were showed in \cite{himho13,holm14,lvw14,wlc16}.
To obtain the desired result, we need to extend the estimate of $\|fg\|_{H^{r-1}(\mathbb{R})}$ for $0\leq r\leq1$ in \cite{himho13}, commonly used in the previous works, to that of $\|fg\|_{H^{r-k}(\mathbb{R})}$ for $0\leq r\leq k$ and $k>1$, which plays a key role in proving the main result.

The rest of the paper is organized as follows. In Section 2, the local well-posedness for (1.1)
with $\sigma=2$ and initial data in $H^{s}(\mathbb{R})\times H^{s-2}(\mathbb{R}), s>\frac{7}{2}$, is established, an explicit lower bound for the maximal existence
time $T$ and an estimate of the solution size are provided.
The H\"{o}lder continuity of the data-to-solution map is showed in Section 3.

Throughout the paper, we denote by $\|\cdot\|_{X}$ the norm of Banach space $X$, $(\cdot,\cdot)$ the inner product of Hilbert space $L^{2}(\mathbb{R})$, and $``\lesssim"$ the inequality up to a positive constant.

\section{Local well-posedness and estimate of the solution size}

In this section, we will give the local well-posedness for Eq.(1.1) with $\sigma=2$, and provide an explicit lower bound for the maximal existence time and an estimate of the solution size.

Setting $\Lambda^{-4}:=(1-\partial_{x}^{2})^{-2}$, the initial-value problem associated to Eq.(1.1) with $\sigma=2$ can be rewritten in the following form:
$$
\begin{array}{l}
\left\{\begin{array}{l}
u_{t}+uu_{x}
 +\partial_{x}\Lambda^{-4}\left(\frac{b}{2}u^{2}+(3-b)u_{x}^{2}-\frac{b+5}{2}u_{xx}^{2}+
 (b-5)u_{x}u_{xxx}+\frac{\kappa}{2}\rho^{2}-\alpha u\right)=0,\quad t>0, ~x\in \mathbb{R},\\[3pt]
\rho_{t}+u\rho_{x}+(b-1)u_{x}\rho=0, \quad t>0, ~x\in \mathbb{R},\\[3pt]
u(0,x)=u_{0}(x), \quad \rho(0,x)=\rho_{0}(x) \quad x\in \mathbb{R},
\end{array}
\right.
\end{array}
\eqno(2.1)
$$

Applying the transport equation theory combined with the method of the Besov spaces, one may obtain the following local well-posedness result for system (2.1), more details can be seen in \cite{zhli17,zhou18}.\\

\noindent\textbf{Theorem 2.1.}  Given $(u_{0}, \rho_{0})\in H^{s}(\mathbb{R})\times H^{s-2}(\mathbb{R}), s>\frac{7}{2}$, there exist a maximal $T=T(u_{0}, \rho_{0})>0$ and a unique solution $(u, \rho)$ to (2.1) such that
$$
(u, \rho)\in C([0,T);H^{s}(\mathbb{R})\times H^{s-2}(\mathbb{R}))\cap C^{1}([0,T);H^{s-1}(\mathbb{R})\times H^{s-3}(\mathbb{R})).
$$
Moreover, the solution depends continuously on the initial data, and $T$ is independent of $s$.\\

Next, we recall the following estimates which will be used later.\\

\noindent\textbf{Lemma 2.1.} (see \cite{kap88}) If $r>0$, then $H^{r}(\mathbb{R})\cap L^{\infty}(\mathbb{R})$
is an algebra. Moreover,
$$
\|fg\|_{H^{r}(\mathbb{R})}
\leq c_{r}(\|f\|_{L^{\infty}(\mathbb{R})}\|g\|_{H^{r}(\mathbb{R})}
 +\|f\|_{H^{r}(\mathbb{R})}\|g\|_{L^{\infty}(\mathbb{R})}),
$$
where $c_{r}$ is a positive constant depending only on $r$.\\

\noindent\textbf{Lemma 2.2.} (see \cite{kap88}) If $r>0$, then
$$
\left\|[\Lambda^{r}, f]g\right\|_{L^{2}(\mathbb{R})}
\leq c_{r}(\|\partial_{x}f\|_{L^{\infty}(\mathbb{R})}\|\Lambda^{r-1}g\|_{L^{2}(\mathbb{R})}
 +\|\Lambda^{r}f\|_{L^{2}(\mathbb{R})}\|g\|_{L^{\infty}(\mathbb{R})}),
$$
where $\Lambda^{r}=(1-\partial_{x}^{2})^{r/2}$ and $c_{r}$ is a positive constant depending only on $r$.\\

\noindent\textbf{Lemma 2.3.} (see \cite{tay14}) If $f\in H^{s}(\mathbb{R})$ with $s>\frac{3}{2}$, then there exists a constant $c>0$ such that for any $g\in L^{2}(\mathbb{R})$ we have
$$
\|[J_{\varepsilon}, f]\partial_{x}g\|_{L^{2}(\mathbb{R})}
\leq c\|f\|_{C^{1}(\mathbb{R})}\|g\|_{L^{2}(\mathbb{R})},
$$
in which for each $\varepsilon\in(0,1]$, the operator $J_{\varepsilon}$ is the Friedrichs mollifier defined by
$$
\begin{array}{l}
J_{\varepsilon}f(x)=j_{\varepsilon}*f(x),
\end{array}
$$
where $j_{\varepsilon}(x)=\frac{1}{\varepsilon}j(\frac{x}{\varepsilon})$ and $j(x)$ is a nonnegative, even, smooth bump function supported in the interval $(-1, 1)$ such that $\int_{\mathbb{R}}j(x)dx=1$.
For any $f\in H^{s}(\mathbb{R})$ with $s\geq0$, we have $J_{\varepsilon}f\rightarrow f$ in $H^{s}(\mathbb{R})$ as $\varepsilon\rightarrow 0$. \\

\noindent\textbf{Theorem 2.2.} Let $(u, \rho)$ be the solution of system (2.1) with initial data $(u_{0}, \rho_{0})\in H^{s}(\mathbb{R})\times H^{s-2}(\mathbb{R}), s> \frac{7}{2}$. Then, the maximal existence time $T$ satisfies
$$
T\geq T_{0}:=\frac{1}{2c_{s}}\ln(1+\frac{1}{\|u_{0}\|_{H^{s}(\mathbb{R})}
+\|\rho_{0}\|_{H^{s-2}(\mathbb{R})}}),
$$
where $c_{s}$ is a constant depending on $s$. Also, we have
$$
\|u\|_{H^{s}(\mathbb{R})}
+\|\rho\|_{H^{s-2}(\mathbb{R})}
\leq 2e^{c_{s}T_{0}}(\|u_{0}\|_{H^{s}(\mathbb{R})}+\|\rho_{0}\|_{H^{s-2}(\mathbb{R})}),
\quad t\in[0, T_{0}].
$$

\noindent\textbf{Proof.}
Note that the products $uu_{x}$, $u\rho_{x}$ only have the regularity of $H^{s-1}(\mathbb{R})$ and $H^{s-3}(\mathbb{R})$ when $(u,\rho)\in H^{s}(\mathbb{R})\times H^{s-2}(\mathbb{R})$. To deal with this problem, we apply the operator $J_{\varepsilon}$ to the system
$$
\begin{array}{rl}
\left\{\begin{array}{l}
(J_{\varepsilon}u)_{t}+J_{\varepsilon}(uu_{x})
+\partial_{x}\Lambda^{-4}[
\frac{b}{2}J_{\varepsilon}(u^{2})+(3-b)J_{\varepsilon}(u_{x}^{2})\\[3pt]
\qquad\qquad-\frac{b+5}{2}J_{\varepsilon}(u_{xx}^{2})+(b-5)J_{\varepsilon}(u_{x}u_{xxx})
+\frac{\kappa}{2}J_{\varepsilon}(\rho^{2})-\alpha J_{\varepsilon}u]=0,\\[5pt]
(J_{\varepsilon}\rho)_{t}+J_{\varepsilon}(u\rho_{x})
+(b-1)J_{\varepsilon}(\rho u_{x})=0.
\end{array}
\right.
\end{array}
\eqno(2.2)
$$
Applying the operator $\Lambda^{s}=(1-\partial_{x}^{2})^{s/2}$ to the first equation of (2.2), then multiplying the resulting equation by $\Lambda^{s}J_{\varepsilon}u$ and integrating with respect to $x\in \mathbb{R}$, we obtain
$$
\begin{array}{rl}
&\frac{1}{2}\frac{d}{dt}\|J_{\varepsilon}u\|_{H^{s}(\mathbb{R})}^{2}
=-\left(\Lambda^{s}J_{\varepsilon}(uu_{x}), \Lambda^{s}J_{\varepsilon}u\right)\\[3pt]
& -\left(\Lambda^{s}J_{\varepsilon}u,\partial_{x}\Lambda^{s}\Lambda^{-4}[\frac{b}{2}J_{\varepsilon}(u^{2})+(3-b)J_{\varepsilon}(u_{x}^{2})
-\frac{b+5}{2}J_{\varepsilon}(u_{xx}^{2})+(b-5)J_{\varepsilon}(u_{x}u_{xxx})
+\frac{\kappa}{2}J_{\varepsilon}(\rho^{2})-\alpha J_{\varepsilon}u]\right).
\end{array}
\eqno(2.3)
$$
In what follows next we use the fact that $\Lambda^{s}$ and $J_{\varepsilon}$ commute and that
$J_{\varepsilon}$ satisfies the properties
$$
(J_{\varepsilon}f, g)=(f, J_{\varepsilon}g)
\quad \mbox{and}\quad
\|J_{\varepsilon}u\|_{H^{s}(\mathbb{R})}\leq \|u\|_{H^{s}(\mathbb{R})}.
$$
Let us estimate the first term of the right hand side of (2.3).
$$
\begin{array}{rl}
&\left|\left(\Lambda^{s}J_{\varepsilon}(uu_{x}), \Lambda^{s}J_{\varepsilon}u\right)\right|\\[3pt]
&=\left|(\Lambda^{s}(uu_{x}),J_{\varepsilon}\Lambda^{s}J_{\varepsilon}u)\right|\\[3pt]
&=\left|([\Lambda^{s}, u]u_{x},J_{\varepsilon}\Lambda^{s}J_{\varepsilon}u)
  +( u\Lambda^{s}u_{x},J_{\varepsilon}\Lambda^{s}J_{\varepsilon}u)\right|\\[3pt]
&=\left|([\Lambda^{s}, u]u_{x},J_{\varepsilon}\Lambda^{s}J_{\varepsilon}u)
  +( J_{\varepsilon}u\partial_{x}\Lambda^{s}u,\Lambda^{s}J_{\varepsilon}u)\right|\\[3pt]
&=\left|([\Lambda^{s}, u]u_{x},J_{\varepsilon}\Lambda^{s}J_{\varepsilon}u)
  +([J_{\varepsilon},u]\partial_{x}\Lambda^{s}u,\Lambda^{s}J_{\varepsilon}u)
  +(uJ_{\varepsilon}\partial_{x}\Lambda^{s}u,\Lambda^{s}J_{\varepsilon}u)\right|\\[3pt]
&\leq\|[\Lambda^{s}, u]u_{x}\|_{L^{2}(\mathbb{R})}
  \|J_{\varepsilon}\Lambda^{s}J_{\varepsilon}u\|_{L^{2}(\mathbb{R})}
  +\|[J_{\varepsilon},u]\partial_{x}\Lambda^{s}u\|_{L^{2}(\mathbb{R})}
  \|\Lambda^{s}J_{\varepsilon}u\|_{L^{2}(\mathbb{R})}\\[3pt]
&\quad+\frac{1}{2}\left|(u_{x}\Lambda^{s}J_{\varepsilon}u,\Lambda^{s}J_{\varepsilon}u)\right|\\[3pt]
&\lesssim \|u\|_{H^{s}(\mathbb{R})}^{3},
\end{array}
$$
where we have used Lemma 2.2 with $r=s$ and Lemma 2.3. Furthermore, we estimate the second term of the right hand side of $(2.3)$ in the following way
$$
\begin{array}{rl}
&\left|\left(\Lambda^{s}J_{\varepsilon}u,\partial_{x}\Lambda^{s}\Lambda^{-4}\left(\frac{b}{2}J_{\varepsilon}(u^{2})+(3-b)J_{\varepsilon}(u_{x}^{2})
-\frac{b+5}{2}J_{\varepsilon}(u_{xx}^{2})+(b-5)J_{\varepsilon}(u_{x}u_{xxx})
+\frac{\kappa}{2}J_{\varepsilon}(\rho^{2})-\alpha J_{\varepsilon}u\right)\right)\right|\\[3pt]
&\leq\left\|\partial_{x}\Lambda^{-4}\left(\frac{b}{2}J_{\varepsilon}(u^{2})+(3-b)J_{\varepsilon}(u_{x}^{2})
-\frac{b+5}{2}J_{\varepsilon}(u_{xx}^{2})+(b-5)J_{\varepsilon}(u_{x}u_{xxx})
+\frac{\kappa}{2}J_{\varepsilon}(\rho^{2})-\alpha J_{\varepsilon}u\right)\right\|_{H^{s}(\mathbb{R})}\|u\|_{H^{s}(\mathbb{R})}\\[3pt]
&\leq\|\frac{b}{2}J_{\varepsilon}(u^{2})+(3-b)J_{\varepsilon}(u_{x}^{2})
-\frac{b+5}{2}J_{\varepsilon}(u_{xx}^{2})+(b-5)J_{\varepsilon}(u_{x}u_{xxx})
+\frac{\kappa}{2}J_{\varepsilon}(\rho^{2})-\alpha J_{\varepsilon}u\|_{H^{s-3}(\mathbb{R})}\|u\|_{H^{s}(\mathbb{R})}\\[3pt]
&\lesssim (\|u\|_{H^{s-3}(\mathbb{R})}^{2}
  +\|u_{x}\|_{H^{s-3}(\mathbb{R})}^{2}
  +\|u_{xx}\|_{H^{s-3}(\mathbb{R})}^{2}
  +\|u_{x}\|_{H^{s-3}(\mathbb{R})}\|u_{xxx}\|_{H^{s-3}(\mathbb{R})}
  +\|\rho\|_{H^{s-3}(\mathbb{R})}^{2}
  +\|u\|_{H^{s-3}(\mathbb{R})})\|u\|_{H^{s}(\mathbb{R})}\\[3pt]
&\lesssim (\|u\|_{H^{s}(\mathbb{R})}^{2}
  +\|\rho\|_{H^{s-2}(\mathbb{R})}^{2}
  +\|u\|_{H^{s}(\mathbb{R})})\|u\|_{H^{s}(\mathbb{R})},
\end{array}
$$
where we have used Lemma 2.1 with $r=s-3$. Thus, we have
$$
\begin{array}{rl}
\frac{1}{2}\frac{d}{dt}\|J_{\varepsilon}u\|_{H^{s}(\mathbb{R})}^{2}
\lesssim (\|u\|_{H^{s}(\mathbb{R})}^{2}
  +\|\rho\|_{H^{s-2}(\mathbb{R})}^{2}
  +\|u\|_{H^{s}(\mathbb{R})})\|u\|_{H^{s}(\mathbb{R})}.
\end{array}
$$
Letting $\varepsilon\rightarrow 0$, we get
$$
\begin{array}{rl}
\frac{d}{dt}\|u\|_{H^{s}(\mathbb{R})}
\lesssim \|u\|_{H^{s}(\mathbb{R})}^{2}
  +\|\rho\|_{H^{s-2}(\mathbb{R})}^{2}
  +\|u\|_{H^{s}(\mathbb{R})}.
\end{array}
\eqno(2.4)
$$

Applying the operator $\Lambda^{s-2}=(1-\partial_{x}^{2})^{(s-2)/2}$ to the second equation of (2.3), then multiplying the resulting equation by $\Lambda^{s-2}J_{\varepsilon}\rho$ and integrating with respect to $x\in \mathbb{R}$, we obtain
$$
\begin{array}{rl}
&\frac{1}{2}\frac{d}{dt}\|J_{\varepsilon}\rho\|_{H^{s-2}(\mathbb{R})}^{2}\\[3pt]
&=-\left(\Lambda^{s-2}J_{\varepsilon}(u\rho_{x}), \Lambda^{s-2}J_{\varepsilon}\rho\right)
-(b-1)\left(\Lambda^{s-2}J_{\varepsilon}(\rho u_{x}), \Lambda^{s-2}J_{\varepsilon}\rho\right)\\[3pt]
&=-\left(\Lambda^{s-2}(u\rho_{x}), J_{\varepsilon}\Lambda^{s-2}J_{\varepsilon}\rho\right)
-(b-1)\left(\Lambda^{s-2}(\rho u_{x}), J_{\varepsilon}\Lambda^{s-2}J_{\varepsilon}\rho\right)\\[3pt]
&=-([\Lambda^{s-2}, u]\rho_{x},J_{\varepsilon}\Lambda^{s-2}J_{\varepsilon}\rho)
  -([J_{\varepsilon},u]\Lambda^{s-2}\rho_{x},\Lambda^{s-2}J_{\varepsilon}\rho)
  -(uJ_{\varepsilon}\Lambda^{s-2}\rho_{x},\Lambda^{s-2}J_{\varepsilon}\rho)\\[3pt]
&\quad
-(b-1)([\Lambda^{s-2}, \rho]u_{x},J_{\varepsilon}\Lambda^{s-2}J_{\varepsilon}\rho)
-(b-1)([J_{\varepsilon},\rho]\Lambda^{s-2}u_{x},\Lambda^{s-2}J_{\varepsilon}\rho)\\[3pt]
&\quad-(b-1)(\rho J_{\varepsilon}\Lambda^{s-2}u_{x},\Lambda^{s-2}J_{\varepsilon}\rho)\\[3pt]
&\lesssim \|[\Lambda^{s-2}, u]\rho_{x}\|_{L^{2}(\mathbb{R})}
\|J_{\varepsilon}\Lambda^{s-2}J_{\varepsilon}\rho\|_{L^{2}(\mathbb{R})}
+\|[J_{\varepsilon},u]\Lambda^{s-2}\rho_{x}\|_{L^{2}(\mathbb{R})}
\|\Lambda^{s-2}J_{\varepsilon}\rho\|_{L^{2}(\mathbb{R})}\\[3pt]
&\quad+|(u_{x}\Lambda^{s-2}J_{\varepsilon}\rho,\Lambda^{s-2}J_{\varepsilon}\rho)|
+\|[\Lambda^{s-2}, \rho]u_{x}\|_{L^{2}(\mathbb{R})}
\|J_{\varepsilon}\Lambda^{s-2}J_{\varepsilon}\rho\|_{L^{2}(\mathbb{R})}\\[3pt]
&\quad+\|[J_{\varepsilon},\rho]\Lambda^{s-2}u_{x}\|_{L^{2}(\mathbb{R})}
\|\Lambda^{s-2}J_{\varepsilon}\rho\|_{L^{2}(\mathbb{R})}
+\|\rho J_{\varepsilon}\Lambda^{s-2}u_{x}\|_{L^{2}(\mathbb{R})}
\|\Lambda^{s-2}J_{\varepsilon}\rho\|_{L^{2}(\mathbb{R})}\\[3pt]
&\lesssim \|u\|_{H^{s}(\mathbb{R})}\|\rho\|_{H^{s-2}(\mathbb{R})}^{2},
\end{array}
$$
where we have used Lemmas 2.2-2.3 and integrating by parts. Letting $\varepsilon\rightarrow 0$, we get
$$
\begin{array}{rl}
\frac{d}{dt}\|\rho\|_{H^{s-2}(\mathbb{R})}
\lesssim \|\rho\|_{H^{s-2}(\mathbb{R})}\|u\|_{H^{s}(\mathbb{R})}.
\end{array}
\eqno(2.5)
$$

Combining (2.4) and (2.5), we have
$$
\begin{array}{rl}
&\frac{d}{dt}(\|u\|_{H^{s}(\mathbb{R})}
+\|\rho\|_{H^{s-2}(\mathbb{R})})\\[3pt]
&\lesssim \|u\|_{H^{s}(\mathbb{R})}^{2}
  +\|\rho\|_{H^{s-2}(\mathbb{R})}^{2}
  +\|\rho\|_{H^{s-2}(\mathbb{R})}\|u\|_{H^{s}(\mathbb{R})}
  +\|u\|_{H^{s}(\mathbb{R})}\\[3pt]
&\leq (\|u\|_{H^{s}(\mathbb{R})}
+\|\rho\|_{H^{s-2}(\mathbb{R})})^{2}+\|u\|_{H^{s}(\mathbb{R})}
+\|\rho\|_{H^{s-2}(\mathbb{R})}.
\end{array}
$$
Letting $y(t)=\|u\|_{H^{s}(\mathbb{R})}
+\|\rho\|_{H^{s-2}(\mathbb{R})}$, then we get
$$
-\frac{d(y^{-1}+1)}{dt}
\leq c_{s}(y^{-1}+1),\quad y_{0}:=y(0)=\|u_{0}\|_{H^{s}(\mathbb{R})}
+\|\rho_{0}\|_{H^{s-2}(\mathbb{R})},
$$
which implies that
$$
y\leq \frac{1}{e^{-c_{s}t}(y_{0}^{-1}+1)-1}.
$$
Setting
$$
T_{0}:=\frac{1}{2c_{s}}\ln(1+\frac{1}{\|u_{0}\|_{H^{s}(\mathbb{R})}
+\|\rho_{0}\|_{H^{s-2}(\mathbb{R})}}),
$$
we see from the above inequality that the solution $(u,\rho)$ exists for $0\leq t\leq T_{0}$ and satisfies a solution size bound
$$
\|u\|_{H^{s}(\mathbb{R})}
+\|\rho\|_{H^{s-2}(\mathbb{R})}
\leq 2e^{c_{s}T_{0}}(\|u_{0}\|_{H^{s}(\mathbb{R})}+\|\rho_{0}\|_{H^{s-2}(\mathbb{R})}),\quad\forall~0\leq t\leq T_{0},
$$
which completes the proof of the theorem. \hfill $\Box$\\

\section{H\"{o}lder continuity}

In this section, we will show that the solution map for system (2.1) is H\"{o}lder continuous
in $H^{s}(\mathbb{R})\times H^{s-2}(\mathbb{R})$, $s>\frac{7}{2}$, equipped with the $H^{r}(\mathbb{R})\times H^{r-2}(\mathbb{R})$-topology for $1\leq r<s$. Firstly, we recall the following lemmas.\\

\noindent\textbf{Lemma 3.1.} (see \cite{taylor03}) If $s>\frac{3}{2}$ and $0\leq \sigma+1\leq s$, then there exists a constant $c>0$ such that
$$
\|[\Lambda^{\sigma}\partial_{x}, f]v\|_{L^{2}(\mathbb{R})}
\leq c\|f\|_{H^{s}(\mathbb{R})}\|v\|_{H^{\sigma}(\mathbb{R})}.
$$

\noindent\textbf{Lemma 3.2.} (see \cite{himkm10}) If $r>\frac{1}{2}$, then there exists a constant $c_{r}>0$ depending only on $r$ such that
$$
\|fg\|_{H^{r-1}(\mathbb{R})}\leq c_{r}\|f\|_{H^{r}(\mathbb{R})}\|g\|_{H^{r-1}(\mathbb{R})}.
$$

Lemma 3.2 gives the estimate of $\|fg\|_{H^{s}(\mathbb{R})}$ for $s>-\frac{1}{2}$, the other cases are provided in the following lemma.\\

\noindent\textbf{Lemma 3.3.} If $0\leq r\leq k$, $j>\frac{1}{2}$ and $j\geq k-r$ with $k\in \mathbb{Z}_{+}$, then
there exists a constant $c_{r,j,k}>0$ depending on $r$, $j$ and $k$ such that
$$
\|fg\|_{H^{r-k}(\mathbb{R})}\leq c_{r,j,k}\|f\|_{H^{j}(\mathbb{R})}\|g\|_{H^{r-k}(\mathbb{R})}.
$$

\noindent\textbf{Proof.} The proof can be done by adapting analogous methods as in \cite{himho13}, in which they only considered the case $k=1$. For the reader's convenience, we provide the arguments with obvious modifications. Similar as the proof of Lemma 3 on $\mathbb{R}$ in \cite{himho13}, we can obtain
$$
\begin{array}{rl}
\|fg\|_{H^{r-k}(\mathbb{R})}^{2}
&=\int_{\mathbb{R}}(1+\xi^{2})^{r-k}|\int_{\mathbb{R}}\widehat{f}(\eta)\widehat{g}(\xi-\eta)d\eta|^{2}d\xi\\[3pt]
&=\int_{\mathbb{R}}(1+\xi^{2})^{r-k}|\int_{\mathbb{R}}(1+\eta^{2})^{\frac{j}{2}}\widehat{f}(\eta)\cdot(1+\eta^{2})^{-\frac{j}{2}}\widehat{g}(\xi-\eta)d\eta|^{2}d\xi\\[3pt]
&\leq\|f\|_{H^{j}(\mathbb{R})}^{2}\int_{\mathbb{R}}|\widehat{g}(\eta)|^{2}
\int_{\mathbb{R}}(1+\xi^{2})^{r-k}(1+(\xi-\eta)^{2})^{-j}d\xi d\eta,
\end{array}
$$
in which we have applied the Cauchy-Schwartz inequality in $\eta$, a change of variables, and changed the order of summation.
To get the desired result, it is sufficient to show that there exists a constant
$c_{r,j,k}>0$ such that
$$
\begin{array}{rl}
\int_{\mathbb{R}}(1+\xi^{2})^{r-k}(1+(\xi-\eta)^{2})^{-j}d\xi
\leq c_{r,j,k}(1+\eta^{2})^{r-k}.
\end{array}
$$
In fact, we can check the inequality under the conditions $j>\frac{1}{2}$ and $j\geq k-r$,
the main difference with proof of Lemma 5 in \cite{himho13} is replacing the discussions on $\frac{1}{2}<r\leq 1$ ($0\leq r<\frac{1}{2}$, $r=\frac{1}{2}$, respectively) by $k-\frac{1}{2}<r\leq k$ ($0\leq r< k-\frac{1}{2}$, $r=k-\frac{1}{2}$, respectively) for the cases $\xi\in [\frac{\eta}{2}, \eta]$ and $\xi\in [\eta,\frac{3\eta}{2}]$.\hfill $\Box$\\

\noindent\textbf{Remark 3.1.} Lemma 3.3 is more general than $(iii)$ of Proposition 2.4 in \cite{guil10} when considering the Sobolev norm of $fg$ with negative index, since it covers the case $j=k-r$ here.\\

\noindent\textbf{Theorem 3.1.} Assume $s>\frac{7}{2}$ and $1\leq r<s$. Then the solution map
for system (2.1) is H\"{o}lder continuous with exponent
$$
\begin{array}{rl}
\beta=
\left\{\begin{array}{l}
1, \quad \mbox{if}~1\leq r\leq s-1~\mbox{and}~s+r\geq 5, \\[3pt]
\frac{2s-5}{s-r}, \quad \mbox{if}~\frac{7}{2}<s<4~\mbox{and}~1\leq r\leq 5-s, \\[3pt]
s-r, \quad \mbox{if}~s-1< r< s
\end{array}
\right.
\end{array}
$$
as a map from $B(0,h):=\{(u, \rho)\in H^{s}(\mathbb{R})\times H^{s-2}(\mathbb{R}): \|u\|_{H^{s}(\mathbb{R})}+\|\rho\|_{H^{s-2}(\mathbb{R})}\leq h\}$ with $H^{r}(\mathbb{R})\times H^{r-2}(\mathbb{R})$-norm to $C([0,T_{0}]; H^{r}(\mathbb{R})\times H^{r-2}(\mathbb{R}))$,
where $T_{0}>0$ is defined as in Theorem 2.2. More precisely, we have
$$
\|(u(t),\rho(t))-(v(t),\theta(t))\|_{C([0,T_{0}]; H^{r}(\mathbb{R})\times H^{r-2}(\mathbb{R}))}
\leq c\|(u(0),\rho(0))-(v(0),\theta(0))\|_{H^{r}(\mathbb{R})\times H^{r-2}(\mathbb{R})}^{\beta},
$$
for all $(u(0),\rho(0)), (v(0),\theta(0))\in B(0,h)$
and $(u(t),\rho(t)), (v(t),\theta(t))$ the solutions corresponding to the initial data $(u(0),\rho(0)), (v(0),\theta(0))$, respectively. The constant $c$ depends on $s, r, T_{0}$ and $h$. \\

\noindent\textbf{Proof.} Define $w=u-v$ and $\eta=\rho-\theta$, then $(w,\eta)$ satisfies that
$$
\begin{array}{rl}
\left\{\begin{array}{l}
w_{t}+\partial_{x}\left(\frac{1}{2}w(u+v)\right)
+\partial_{x}\Lambda^{-4}[\frac{b}{2}w(u+v)+(3-b)w_{x}(u_{x}+v_{x})
-\frac{b+5}{2}w_{xx}(u_{xx}+v_{xx})\\[5pt]
\qquad +(b-5)w_{x}u_{xxx}+(b-5)v_{x}w_{xxx}+\frac{\kappa}{2}\eta(\rho+\theta)-\alpha w]=0, \quad t>0, ~x\in \mathbb{R},\\[5pt]
\eta_{t}+u\eta_{x}+w\theta_{x}=-(b-1)(w_{x}\rho+v_{x}\eta), \quad t>0, ~x\in \mathbb{R},\\[5pt]
w(0,x)=u_{0}-v_{0},~\eta(0,x)=\rho_{0}-\theta_{0}, \quad x\in \mathbb{R}.
\end{array}
\right.
\end{array}
\eqno(3.1)
$$

$(i)$ We first consider the case $1\leq r\leq s-1$ and $r+s\geq 5$, where $s>\frac{7}{2}$. Applying $\Lambda^{r}$ to the first equation of (3.1), then multiplying both sides by $\Lambda^{r}w$ and integrating over $\mathbb{R}$ with respect to $x$, we get
$$
\begin{array}{rl}
\frac{1}{2}\frac{d}{dt}\|w\|_{H^{r}(\mathbb{R})}^{2}
&=-\left(\Lambda^{r}\partial_{x}
  (\frac{1}{2}w(u+v)),\Lambda^{r}w\right)\\[3pt]
&\quad-(\Lambda^{r}\partial_{x}\Lambda^{-4}[\frac{b}{2}w(u+v)+(3-b)w_{x}(u_{x}+v_{x})
-\frac{b+5}{2}w_{xx}(u_{xx}+v_{xx})\\[3pt]
&\qquad+(b-5)w_{x}u_{xxx}+(b-5)v_{x}w_{xxx}+\frac{k}{2}\eta(\rho+\theta)-\alpha w], \Lambda^{r}w)\\[3pt]
&:=E_{1}+E_{2}.
\end{array}
$$
To get the desired result, we need to estimate $E_{1}$ and $E_{2}$.

\emph{Estimate} $E_{1}$. By using Lemma 3.1, integrating by parts and the Sobolev embedding theorem $H^{r}(\mathbb{R})\hookrightarrow L^{\infty}(\mathbb{R})$ for $r>\frac{1}{2}$, we have
$$
\begin{array}{rl}
|E_{1}|
&=\left|-\left(\Lambda^{r}\partial_{x}
  (\frac{1}{2}w(u+v)),\Lambda^{r}wdx\right)\right|\\[5pt]
&=\left|-\left([\Lambda^{r}\partial_{x}, \frac{1}{2}(u+v)]w, \Lambda^{r}w\right)
  -\left(\frac{1}{2}(u+v)\Lambda^{r}\partial_{x}w,\Lambda^{r}w\right)\right|\\[5pt]
&=\left|-\left([\Lambda^{r}\partial_{x}, \frac{1}{2}(u+v)]w, \Lambda^{r}w\right)
  +\left(\frac{1}{4}\partial_{x}(u+v)\Lambda^{r}w,\Lambda^{r}w\right)\right|\\[5pt]
&\lesssim\|u+v\|_{H^{s}(\mathbb{R})}\|w\|_{H^{r}(\mathbb{R})}^{2}.
\end{array}
$$

\emph{Estimate} $E_{2}$. It is easy to show that
$$
\begin{array}{rl}
|E_{2}|
&=|-(\Lambda^{r}\partial_{x}\Lambda^{-4}[\frac{b}{2}w(u+v)+(3-b)w_{x}(u_{x}+v_{x})
-\frac{b+5}{2}w_{xx}(u_{xx}+v_{xx})\\[5pt]
&\qquad+(b-5)w_{x}u_{xxx}+(b-5)v_{x}w_{xxx}+\frac{\kappa}{2}\eta(\rho+\theta)-\alpha w], \Lambda^{r}w)|\\[5pt]
&\leq \|\partial_{x}\Lambda^{-4}[\frac{b}{2}w(u+v)+(3-b)w_{x}(u_{x}+v_{x})
-\frac{b+5}{2}w_{xx}(u_{xx}+v_{xx})\\[5pt]
&\qquad+(b-5)w_{x}u_{xxx}+(b-5)v_{x}w_{xxx}+\frac{\kappa}{2}\eta(\rho+\theta)-\alpha w]\|_{H^{r}(\mathbb{R})}\|w\|_{H^{r}(\mathbb{R})}.
\end{array}
$$
Using integrating by parts, we have
$$
\begin{array}{rl}
&\|\partial_{x}\Lambda^{-4}[\frac{b}{2}w(u+v)+(3-b)w_{x}(u_{x}+v_{x})
-\frac{b+5}{2}w_{xx}(u_{xx}+v_{xx})\\[5pt]
&\qquad+(b-5)w_{x}u_{xxx}+(b-5)v_{x}w_{xxx}+\frac{\kappa}{2}\eta(\rho+\theta)-\alpha w]\|_{H^{r}(\mathbb{R})}\\[5pt]
&=\|\partial_{x}\Lambda^{-4}[\frac{b}{2}w(u+v)+(3-b)w_{x}(u_{x}+v_{x})
+w_{xx}(\frac{5-3b}{2}u_{xx}-\frac{b+5}{2}v_{xx})\\[5pt]
&\qquad
+(b-5)v_{x}w_{xxx}+\frac{\kappa}{2}\eta(\rho+\theta)-\alpha w]
+(b-5)\partial_{x}^{2}\Lambda^{-4}(w_{x}u_{xx})\|_{H^{r}(\mathbb{R})}\\[5pt]
&\lesssim [\|w(u+v)\|_{H^{r-3}(\mathbb{R})}
+\|w_{x}(u_{x}+v_{x})\|_{H^{r-3}(\mathbb{R})}
+\|w_{xx}(u_{xx}+v_{xx})\|_{H^{r-3}(\mathbb{R})}\\[5pt]
&\qquad+\|v_{x}w_{xxx}\|_{H^{r-3}(\mathbb{R})}
+\|\eta(\rho+\theta)\|_{H^{r-3}(\mathbb{R})}
+\|w\|_{H^{r-3}(\mathbb{R})}]
+\|w_{x}u_{xx}\|_{H^{r-2}(\mathbb{R})}\\[5pt]
&:=F_{1}+F_{2}.
\end{array}
$$

For $F_{1}$, if $r>\frac{5}{2}$, we have
$$
\begin{array}{rl}
F_{1}
&\lesssim
\|w\|_{H^{r-3}(\mathbb{R})}\|u+v\|_{H^{r-2}(\mathbb{R})}
+\|w_{x}\|_{H^{r-3}(\mathbb{R})}\|u_{x}+v_{x}\|_{H^{r-2}(\mathbb{R})}\\[5pt]
&\quad+\|w_{xx}\|_{H^{r-3}(\mathbb{R})}\|u_{xx}+v_{xx}\|_{H^{r-2}(\mathbb{R})}
+\|w_{xxx}\|_{H^{r-3}(\mathbb{R})}\|v_{x}\|_{H^{r-2}(\mathbb{R})}\\[5pt]
&\quad+\|\eta\|_{H^{r-3}(\mathbb{R})}\|\rho+\theta\|_{H^{r-2}(\mathbb{R})}
+\|w\|_{H^{r-3}(\mathbb{R})}\\[5pt]
&\lesssim \|w\|_{H^{r}(\mathbb{R})}(\|u\|_{H^{s}(\mathbb{R})}+\|v\|_{H^{s}(\mathbb{R})}+1)
+\|\eta\|_{H^{r-2}(\mathbb{R})}(\|\rho\|_{H^{s-2}(\mathbb{R})}+\|\theta\|_{H^{s-2}(\mathbb{R})})
\end{array}
$$
by using Lemma 3.2 and the fact $r\leq s-1$.

For $1\leq r\leq \frac{5}{2}$, applying Lemma 3.3 with $k=3$ to the term $F_{1}$,
we have
$$
\begin{array}{rl}
F_{1}
&\lesssim
\|w\|_{H^{r-3}(\mathbb{R})}\|u+v\|_{H^{j}(\mathbb{R})}
+\|w_{x}\|_{H^{r-3}(\mathbb{R})}\|u_{x}+v_{x}\|_{H^{j}(\mathbb{R})}
+\|w_{xx}\|_{H^{r-3}(\mathbb{R})}\|u_{xx}+v_{xx}\|_{H^{j}(\mathbb{R})}\\[5pt]
&\quad+\|w_{xxx}\|_{H^{r-3}(\mathbb{R})}\|v_{x}\|_{H^{j}(\mathbb{R})}
+\|\eta\|_{H^{r-3}(\mathbb{R})}\|\rho+\theta\|_{H^{j}(\mathbb{R})}
+\|w\|_{H^{r-3}(\mathbb{R})}\\[5pt]
&\leq\|w\|_{H^{r-3}(\mathbb{R})}(\|u\|_{H^{j}(\mathbb{R})}+\|v\|_{H^{j}(\mathbb{R})})
+\|w\|_{H^{r-2}(\mathbb{R})}(\|u\|_{H^{j+1}(\mathbb{R})}+\|v\|_{H^{j+1}(\mathbb{R})})\\[3pt]
&\quad+\|w\|_{H^{r-1}(\mathbb{R})}(\|u\|_{H^{j+2}(\mathbb{R})}+\|v\|_{H^{j+2}(\mathbb{R})})
+\|w\|_{H^{r}(\mathbb{R})}\|v\|_{H^{j+1}(\mathbb{R})}\\[5pt]
&\quad+\|\eta\|_{H^{r-3}(\mathbb{R})}(\|\rho\|_{H^{j}(\mathbb{R})}+\|\theta\|_{H^{j}(\mathbb{R})})
+\|w\|_{H^{r-3}(\mathbb{R})},
\end{array}
$$
where $j$ satisfies $j>\frac{1}{2}$ and $j\geq3-r$. Since $s\geq 5-r$, we can take
$j=s-2$, and then
$$
\begin{array}{rl}
F_{1}\lesssim \|w\|_{H^{r}(\mathbb{R})}(\|u\|_{H^{s}(\mathbb{R})}+\|v\|_{H^{s}(\mathbb{R})}+1)
+\|\eta\|_{H^{r-2}(\mathbb{R})}(\|\rho\|_{H^{s-2}(\mathbb{R})}+\|\theta\|_{H^{s-2}(\mathbb{R})}).
\end{array}
$$

For $F_{2}$, if $r>\frac{3}{2}$, we have
$$
F_{2}
\lesssim\|w_{x}\|_{H^{r-2}(\mathbb{R})}\|u_{xx}\|_{H^{r-1}(\mathbb{R})}
\leq \|w\|_{H^{r}(\mathbb{R})}\|u\|_{H^{s}(\mathbb{R})},
$$
by using Lemma 3.2 and the fact $r\leq s-1$.

For $1\leq r\leq \frac{3}{2}$, applying Lemma 3.3 with $k=2$ to the term $F_{2}$,
we have
$$
\begin{array}{rl}
F_{2}
\lesssim \|w_{x}\|_{H^{r-2}(\mathbb{R)}}\|u_{xx}\|_{H^{l}(\mathbb{R)}}
\leq \|w\|_{H^{r-1}(\mathbb{R})}\|u\|_{H^{l+2}(\mathbb{R})},
\end{array}
\eqno(3.2)
$$
where $l$ satisfies $l>\frac{1}{2}$ and $l\geq 2-r$. Since $s\geq 5-r$, we can take
$l=s-3$, and then
$$
F_{2}
\lesssim \|w\|_{H^{r}(\mathbb{R})}\|u\|_{H^{s}(\mathbb{R})}.
$$
Thus,
$$
\begin{array}{rl}
\frac{1}{2}\frac{d}{dt}\|w\|_{H^{r}(\mathbb{R})}^{2}
&\lesssim \|w\|_{H^{r}(\mathbb{R})}^{2}
(\|u\|_{H^{s}(\mathbb{R})}+\|v\|_{H^{s}(\mathbb{R})}+1)\\[3pt]
&\quad+\|w\|_{H^{r}(\mathbb{R})}\|\eta\|_{H^{r-2}(\mathbb{R})}
(\|\rho\|_{H^{s-2}(\mathbb{R})}+\|\theta\|_{H^{s-2}(\mathbb{R})}).
\end{array}
\eqno(3.3)
$$

On the other hand, applying $\Lambda^{r-2}$ to the second equation of (3.1), then multiplying both sides by $\Lambda^{r-2}\eta$ and integrating over $\mathbb{R}$ with respect to $x$, we get
$$
\begin{array}{rl}
&\frac{1}{2}\frac{d}{dt}\|\eta\|_{H^{r-2}(\mathbb{R})}^{2}\\[3pt]
&=-\left(\Lambda^{r-2}(u\eta_{x}),\Lambda^{r-2}\eta\right)
-\left(\Lambda^{r-2}(w\theta_{x}),\Lambda^{r-2}\eta\right)
-(b-1)\left(\Lambda^{r-2}(w_{x}\rho+v_{x}\eta),\Lambda^{r-2}\eta\right)\\[3pt]
&=-\left(\Lambda^{r-2}\partial_{x}(u\eta),\Lambda^{r-2}\eta\right)
-\left(\Lambda^{r-2}(w\theta_{x}-u_{x}\eta),\Lambda^{r-2}\eta\right)
-(b-1)\left(\Lambda^{r-2}(w_{x}\rho+v_{x}\eta),\Lambda^{r-2}\eta\right)\\[3pt]
&=-\left([\Lambda^{r-2}\partial_{x},u]\eta,\Lambda^{r-2}\eta\right)
-\left(u\Lambda^{r-2}\partial_{x}\eta,\Lambda^{r-2}\eta\right)
-\left(\Lambda^{r-2}(w\theta_{x}-u_{x}\eta),\Lambda^{r-2}\eta\right)\\[3pt]
&\quad-(b-1)\left(\Lambda^{r-2}(w_{x}\rho+v_{x}\eta),\Lambda^{r-2}\eta\right)\\[3pt]
&=-\left([\Lambda^{r-2}\partial_{x},u]\eta,\Lambda^{r-2}\eta\right)
+\frac{1}{2}\left(u_{x}\Lambda^{r-2}\eta,\Lambda^{r-2}\eta\right)
-\left(\Lambda^{r-2}(w\theta_{x}-u_{x}\eta),\Lambda^{r-2}\eta\right)\\[3pt]
&\quad-(b-1)\left(\Lambda^{r-2}(w_{x}\rho+v_{x}\eta),\Lambda^{r-2}\eta\right)\\[3pt]
&\lesssim \left(\|[\Lambda^{r-2}\partial_{x},u]\eta\|_{L^{2}(\mathbb{R})}
+\|u_{x}\Lambda^{r-2}\eta\|_{L^{2}(\mathbb{R})}
+\|w\theta_{x}-u_{x}\eta\|_{H^{r-2}(\mathbb{R})}
+\|w_{x}\rho+v_{x}\eta\|_{H^{r-2}(\mathbb{R})}\right)\|\eta\|_{H^{r-2}(\mathbb{R})}.
\end{array}
$$
Since $r\geq 1$, by Lemma 3.1, we know
$$
\|[\Lambda^{r-2}\partial_{x},u]\eta\|_{L^{2}(\mathbb{R})}
\lesssim \|u\|_{H^{s}(\mathbb{R})}\|\eta\|_{H^{r-2}(\mathbb{R})}.
$$
If $r>\frac{3}{2}$, by Lemma 3.2, we have
$$
\begin{array}{rl}
&\|w\theta_{x}-u_{x}\eta\|_{H^{r-2}(\mathbb{R})}
+\|w_{x}\rho+v_{x}\eta\|_{H^{r-2}(\mathbb{R})}\\[3pt]
&\lesssim \|w\|_{H^{r-1}(\mathbb{R})}\|\theta_{x}\|_{H^{r-2}(\mathbb{R})}
+\|u_{x}\|_{H^{r-1}(\mathbb{R})}\|\eta\|_{H^{r-2}(\mathbb{R})}
+\|w_{x}\|_{H^{r-1}(\mathbb{R})}\|\rho\|_{H^{r-2}(\mathbb{R})}\\[3pt]
&\quad+\|v_{x}\|_{H^{r-1}(\mathbb{R})}\|\eta\|_{H^{r-2}(\mathbb{R})}\\[3pt]
&\leq \|w\|_{H^{r}(\mathbb{R})}(\|\theta\|_{H^{s-2}(\mathbb{R})}+\|\rho\|_{H^{s-2}(\mathbb{R})})
+\|\eta\|_{H^{r-2}(\mathbb{R})}(\|u\|_{H^{s}(\mathbb{R})}+\|v\|_{H^{s}(\mathbb{R})}).
\end{array}
$$
If $1\leq r\leq\frac{3}{2}$, similar to (3.2), we have
$$
\begin{array}{rl}
&\|w\theta_{x}-u_{x}\eta\|_{H^{r-2}(\mathbb{R})}
+\|w_{x}\rho+v_{x}\eta\|_{H^{r-2}(\mathbb{R})}\\[3pt]
&\lesssim \|w\|_{H^{r-2}(\mathbb{R})}\|\theta_{x}\|_{H^{s-3}(\mathbb{R})}
+\|u_{x}\|_{H^{s-3}(\mathbb{R})}\|\eta\|_{H^{r-2}(\mathbb{R})}
+\|w_{x}\|_{H^{r-2}(\mathbb{R})}\|\rho\|_{H^{s-3}(\mathbb{R})}\\[3pt]
&\quad+\|v_{x}\|_{H^{s-3}(\mathbb{R})}\|\eta\|_{H^{r-2}(\mathbb{R})}\\[3pt]
&\leq \|w\|_{H^{r}(\mathbb{R})}(\|\theta\|_{H^{s-2}(\mathbb{R})}+\|\rho\|_{H^{s-2}(\mathbb{R})})
+\|\eta\|_{H^{r-2}(\mathbb{R})}(\|u\|_{H^{s}(\mathbb{R})}+\|v\|_{H^{s}(\mathbb{R})}).
\end{array}
$$
Moreover, it is easy to get
$$
\|u_{x}\Lambda^{r-2}\eta\|_{L^{2}(\mathbb{R})}
\leq \|u_{x}\|_{L^{\infty}(\mathbb{R})}\|\Lambda^{r-2}\eta\|_{L^{2}(\mathbb{R})}
\leq \|u\|_{H^{s}(\mathbb{R})}\|\eta\|_{H^{r-2}(\mathbb{R})}.
$$
Thus, we have
$$
\begin{array}{rl}
&\frac{1}{2}\frac{d}{dt}\|\eta\|_{H^{r-2}(\mathbb{R})}^{2}\\[3pt]
&\lesssim \|w\|_{H^{r}(\mathbb{R})}\|\eta\|_{H^{r-2}(\mathbb{R})}(\|\theta\|_{H^{s-2}(\mathbb{R})}+\|\rho\|_{H^{s-2}(\mathbb{R})})
+\|\eta\|_{H^{r-2}(\mathbb{R})}^{2}(\|u\|_{H^{s}(\mathbb{R})}+\|v\|_{H^{s}(\mathbb{R})}).
\end{array}
\eqno(3.4)
$$

Combing (3.3) and (3.4), and using the solution size estimate in Theorem 2.2, we get
$$
\begin{array}{rl}
&\frac{d}{dt}(\|w\|_{H^{r}(\mathbb{R})}+\|\eta\|_{H^{r-2}(\mathbb{R})})\\[3pt]
&\lesssim (\|w\|_{H^{r}(\mathbb{R})}+\|\eta\|_{H^{r-2}(\mathbb{R})})
(\|u\|_{H^{s}(\mathbb{R})}+\|v\|_{H^{s}(\mathbb{R})}
+\|\theta\|_{H^{s-2}(\mathbb{R})}+\|\rho\|_{H^{s-2}(\mathbb{R})})\\[3pt]
&\lesssim 4e^{c_{s}T_{0}}h(\|w\|_{H^{r}(\mathbb{R})}+\|\eta\|_{H^{r-2}(\mathbb{R})}),
\end{array}
$$
which implies that
$$
\begin{array}{l}
\|w\|_{H^{r}(\mathbb{R})}+\|\eta\|_{H^{r-2}(\mathbb{R})}\leq e^{CT_{0}}(\|w_{0}\|_{H^{r}(\mathbb{R})}+\|\eta_{0}\|_{H^{r-2}(\mathbb{R})}),
\end{array}
\eqno(3.5)
$$
where $C$ is a constant depending on $s, r, T_{0}$ and $h$.

$(ii)$ Next, we consider the case $\frac{7}{2}<s<4$ and $1\leq r\leq 5-s$.
By the condition $r\leq 5-s$ and (3.5), we have
$$
\begin{array}{l}
\|w\|_{H^{r}(\mathbb{R})}+\|\eta\|_{H^{r-2}(\mathbb{R})}\leq
\|w\|_{H^{5-s}(\mathbb{R})}+\|\eta\|_{H^{5-s}(\mathbb{R})}\leq e^{CT_{0}}(\|w_{0}\|_{H^{5-s}(\mathbb{R})}+\|\eta_{0}\|_{H^{3-s}(\mathbb{R})}).
\end{array}
$$
Using interpolation inequalities, we obtain
$$
\begin{array}{rl}
\|w_{0}\|_{H^{5-s}(\mathbb{R})}+\|\eta_{0}\|_{H^{3-s}(\mathbb{R})}
&\leq \|w_{0}\|_{H^{r}(\mathbb{R})}^{\frac{2s-5}{s-r}}
\|w_{0}\|_{H^{s}(\mathbb{R})}^{\frac{5-s-r}{s-r}}
+\|\eta_{0}\|_{H^{r-2}(\mathbb{R})}^{\frac{2s-5}{s-r}}
\|\eta_{0}\|_{H^{s-2}(\mathbb{R})}^{\frac{5-s-r}{s-r}}\\[3pt]
&\leq (\|w_{0}\|_{H^{r}(\mathbb{R})}^{\frac{2s-5}{s-r}}
+\|\eta_{0}\|_{H^{r-2}(\mathbb{R})}^{\frac{2s-5}{s-r}})
(\|w_{0}\|_{H^{s}(\mathbb{R})}^{\frac{5-s-r}{s-r}}
+\|\eta_{0}\|_{H^{s-2}(\mathbb{R})}^{\frac{5-s-r}{s-r}})\\[3pt]
&\leq 4(\|w_{0}\|_{H^{r}(\mathbb{R})}
+\|\eta_{0}\|_{H^{r-2}(\mathbb{R})})^{\frac{2s-5}{s-r}}
(\|w_{0}\|_{H^{s}(\mathbb{R})}
+\|\eta_{0}\|_{H^{s-2}(\mathbb{R})})^{\frac{5-s-r}{s-r}}\\[3pt]
&\leq 4(2h)^{\frac{5-s-r}{s-r}}(\|w_{0}\|_{H^{r}(\mathbb{R})}
+\|\eta_{0}\|_{H^{r-2}(\mathbb{R})})^{\frac{2s-5}{s-r}}.
\end{array}
$$
Thus, we get
$$
\begin{array}{l}
\|w\|_{H^{r}(\mathbb{R})}+\|\eta\|_{H^{r-2}(\mathbb{R})}
\lesssim (\|w_{0}\|_{H^{r}(\mathbb{R})}
+\|\eta_{0}\|_{H^{r-2}(\mathbb{R})})^{\frac{2s-5}{s-r}}.
\end{array}
$$

$(iii)$ Now we consider the case $s-1<r<s$, where $s>\frac{7}{2}$. Using interpolation inequalities, we obtain
$$
\begin{array}{rl}
\|w\|_{H^{r}(\mathbb{R})}+\|\eta\|_{H^{r-2}(\mathbb{R})}
&\leq \|w\|_{H^{s-1}(\mathbb{R})}^{s-r}
\|w\|_{H^{s}(\mathbb{R})}^{1+r-s}
+\|\eta\|_{H^{s-3}(\mathbb{R})}^{s-r}
\|\eta\|_{H^{s-2}(\mathbb{R})}^{1+r-s}\\[3pt]
&\leq (\|w\|_{H^{s-1}(\mathbb{R})}^{s-r}
+\|\eta\|_{H^{s-3}(\mathbb{R})}^{s-r})
(\|w\|_{H^{s}(\mathbb{R})}^{1+r-s}
+\|\eta\|_{H^{s-2}(\mathbb{R})}^{1+r-s})\\[3pt]
&\leq 4(\|w\|_{H^{s-1}(\mathbb{R})}
+\|\eta\|_{H^{s-3}(\mathbb{R})})^{s-r}
(\|w\|_{H^{s}(\mathbb{R})}
+\|\eta\|_{H^{s-2}(\mathbb{R})})^{1+r-s}.
\end{array}
\eqno(3.6)
$$
By applying inequality (3.5), we have
$$
\begin{array}{l}
\|w\|_{H^{s-1}(\mathbb{R})}
+\|\eta\|_{H^{s-3}(\mathbb{R})}
\leq e^{CT_{0}}(\|w_{0}\|_{H^{s-1}(\mathbb{R})}
+\|\eta_{0}\|_{H^{s-3}(\mathbb{R})}).
\end{array}
\eqno(3.7)
$$
Also, using the solution size estimate in Theorem 2.2, we get
$$
\begin{array}{rl}
\|w\|_{H^{s}(\mathbb{R})}
+\|\eta\|_{H^{s-2}(\mathbb{R})}
&\leq (\|u\|_{H^{s}(\mathbb{R})}+\|\rho\|_{H^{s-2}(\mathbb{R})})
+(\|v\|_{H^{s}(\mathbb{R})}+\|\theta\|_{H^{s-2}(\mathbb{R})})\\[3pt]
&\leq 2e^{c_{s}T_{0}}(\|u_{0}\|_{H^{s}(\mathbb{R})}+\|\rho_{0}\|_{H^{s-2}(\mathbb{R})}
+\|v_{0}\|_{H^{s}(\mathbb{R})}+\|\theta_{0}\|_{H^{s-2}(\mathbb{R})})\\[3pt]
&\leq 4e^{cT_{0}}h.
\end{array}
\eqno(3.8)
$$
Combining (3.6), (3.7) and (3.8) gives
$$
\begin{array}{l}
\|w\|_{H^{r}(\mathbb{R})}+\|\eta\|_{H^{r-2}(\mathbb{R})}
\lesssim (\|w_{0}\|_{H^{r}(\mathbb{R})}
+\|\eta_{0}\|_{H^{r-2}(\mathbb{R})})^{s-r}.
\end{array}
$$
This completes the proof of Theorem 3.1.\hfill $\Box$\\

\noindent\textbf{Remark 3.2.} If $\rho\equiv 0$, then the condition $1\leq r<s$ in
Theorem 3.1 can be extended to $0\leq r<s$, and the exponent $\beta$ is defined as follows
$$
\begin{array}{rl}
\beta=
\left\{\begin{array}{l}
1, \quad \mbox{if}~0\leq r\leq s-1~\mbox{and}~s+r\geq 5, \\[3pt]
\frac{2s-5}{s-r}, \quad \mbox{if}~\frac{7}{2}<s<5~\mbox{and}~0\leq r\leq 5-s, \\[3pt]
s-r, \quad \mbox{if}~s-1< r< s.
\end{array}
\right.
\end{array}
$$

\section*{Acknowledgments}
Wang's work is supported by the Fundamental Research Funds for the Central Universities.
Li's work is supported by the NSFC (No:11571057).
\label{}





\bibliographystyle{model3-num-names}
\bibliography{<your-bib-database>}

\begin{thebibliography}{s42}

\bibitem{brc07}
A. Bressan, A. Constantin,
Global dissipative solutions of the Camassa-Holm equation,
Anal. Appl. 5 (2007) 1--27.

\bibitem{ch93}
R. Camassa, D.D. Holm,
An integrable shallow water equation with peaked solitons,
Phys. Rev. Lett. 71 (1993) 1661--1664.

\bibitem{chlz06}
M. Chen, S. Liu, Y. Zhang,
A 2-component generalization of the Camassa-Holm equation and its solutions,
Lett. Math. Phys. 7 (2006) 1--15.

\bibitem{chzh17}
R. Chen, S. Zhou,
Well-posedness and persistence properties for two-component higher order Camassa-Holm systems with fractional inertia operator,
Nonlinear Anal.: RWA 33 (2017) 121--138.

\bibitem{chlz13}
R.M. Chen, Y. Liu, P. Zhang,
The H\"{o}lder continuity of the solution map to the $b$-family equation in weak topology,
Math. Ann. 357 (2013) 1245--1289.

\bibitem{ce98}
A. Constantin, J. Escher,
Wave breaking for nonlinear nonlocal shallow water equations,
Acta Math. 181 (1998) 229--243.

\bibitem{ces98}
A. Constantin, J. Escher,
Global existence and blow-up for a shallow water equation,
Ann. Scuola Norm. Sup. Pisa 26 (1998) 303--328.

\bibitem{cm99}
A. Constantin, H.P. McKean,
A shallow water equation on the circle,
Comm. Pure Appl. Math. 52 (1999) 949--982.

\bibitem{cm00}
A. Constantin, L. Molinet,
Global weak solutions for a shallow water equation,
Comm. Math. Phys. 211 (2000) 45--61.

\bibitem{desp99}
A. Degasperis, M. Procesi,
Asymptotic integrability, Symmetry and Perturbation Theory, World Scientific, Singapore, 1999, pp. 23--37.

\bibitem{ehkl16}
J. Escher, D. Henry, B. Kolev, T. Lyons,
Two-component equations modelling water waves with constant vorticity,
Ann. Mat. Pura Appl. 195 (2016) 249--271.

\bibitem{esko14}
J. Escher, B. Kolev,
Geodesic completeness for Sobolev $H^{s}$-metrics on the diffeomorphism group of the circle,
J. Evol. Equ. 14 (2014) 949--968.

\bibitem{eskol14}
J. Escher, B. Kolev,
Right-invariant Sobolev metrics of fractional order on the diffeomorphism group of the circle,
J. Geom. Mech. 6 (2014) 335--372.

\bibitem{elyi07}
J. Escher, O. Lechtenfeld, Z. Yin,
Well-posedness and blow-up phenomena for the 2-component Camassa-Holm equation,
Discrete Contin. Dyn. Syst. 19 (2007) 493--513.

\bibitem{eslyin06}
J. Escher, Y. Liu, Z. Yin,
Global weak solutions and blow-up structure for the Degasperis-Procesi equation,
J. Funct. Anal. 241 (2006) 457--485.

\bibitem{esly15}
J. Escher, T. Lyons,
Two-component higher order Camassa-Holm systems with fractional inertia operator: A geometric approach, J. Geom. Mech. 7 (2015) 281--293.

\bibitem{esyi08}
J. Escher, Z. Yin,
Well-posedness, blow-up phenomena, and global solutions for the $b$-equation,
J. Reine Angew. Math. 624 (2008) 51--80.

\bibitem{guhey15}
C. Guan, H. He, Z. Yin,
Well-posedness, blow-up phenomena and persistence properties for a two-component water wave system,
Nonlinear Anal.: RWA 25 (2015) 219--237.

\bibitem{guiyin10}
C. Guan, Z. Yin,
Global existence and blow-up phenomena for an integrable two component Camassa-Holm shallow water system,
J. Diff. Equ. 248 (2010) 2003--2014.

\bibitem{guiyin11}
C. Guan, Z. Yin,
Global weak solutions for a two-component Camassa-Holm shallow water system,
J. Funct. Anal. 260 (2011) 1132--1154.

\bibitem{guil10}
G. Gui, Y. Liu,
On the global existence and wave-breaking criteria for the two-component Camassa-Holm system,
J. Funct. Anal. 258 (2010) 4251--4278.

\bibitem{heyin14}
H. He, Z. Yin,
On the Cauchy problem for a generalized two-component shallow water wave system with fractional
higher-order inertia operators,
Discrete Contin. Dyn. Syst. 37 (2017) 1509--1537.

\bibitem{himho13}
A. Himonas, J. Holmes,
H\"{o}lder continuity of the solution map for the Novikov equation,
J. Math. Phys. 54 (2013) 1--11.

\bibitem{hmkz07}
A. Himonas, G. Misio{\l}ek, C. Kenig, Y. Zhou,
Persistence properties and unique continuation of solutions of the Camassa-Holm equation,
Comm. Math. Phys. 271 (2007) 511--522.

\bibitem{himkm10}
A. Himonas, C. Kenig, G. Misio{\l}ek,
Non-uniform dependence for the periodic CH equation,
Comm. Partial Differential Equations 35 (2010) 1145--1162.

\bibitem{holm14}
J. Holmes,
Continuity properties of the data-to-solution map for the generalized Camassa-Holm equation,
J. Math. Anal. Appl. 417 (2014) 635--642.

\bibitem{kap88}
T. Kato, G. Ponce,
Commutator estimates and the Euler and Navier-Stokes equations,
Comm. Pure Appl. Math. 41 (1988) 203--208.

\bibitem{lene05a}
J. Lenells,
Traveling wave solutions of the Camassa-Holm equation,
J. Diff. Equ. 217 (2005) 393--430.

\bibitem{lene05b}
J. Lenells,
Traveling wave solutions of the Degasperis-Procesi equation,
J. Math. Anal. Appl. 306 (2005) 72--82.

\bibitem{lo00}
Y. Li, P. Olver,
Well-posedness and blow-up solutions for an integrable nonlinearly dispersive model wave equation,
J. Diff. Equ. 162 (2000) 27--63.

\bibitem{lvw14}
G. Lv, X. Wang,
H\"{o}lder continuity on $\mu$-b equation,
Nonlinear Anal. 102 (2014) 30--35.

\bibitem{mczh09}
R. McLachlan, X. Zhang,
Well-posedness of a modified Camassa-Holm equations,
J. Diff. Equ. 246 (2009) 3241--3259.

\bibitem{tay14}
M. Taylor,
Pseudodifferential Operators and Nonlinear PDE,
Birkh\"{a}user, Boston (1991)

\bibitem{taylor03}
M. Taylor,
Commutator estimates,
Proc. Amer. Math. Soc. 131 (2003) 1501--1507.

\bibitem{wlc16}
F. Wang, F. Li, Q. Chen,
On the Cauchy problem for a weakly dissipative generalized $\mu$-Hunter-Saxton equation,
Monatsh. Math. 181 (2016) 715--744.

\bibitem{xz00}
Z. Xin, P. Zhang,
On the weak solutions to a shallow water equation,
Comm. Pure Appl. Math. 53 (2000) 1411--1433.

\bibitem{zhli17}
L. Zhang, X. Li,
The local well-posedness, blow-up criteria and Gevrey regularity of solutions for a two-component high-order Camassa-Holm system,
Nonlinear Anal.: RWA 35 (2017) 414--440.

\bibitem{zhou18}
S. Zhou,
Well-posedness, blow-up phenomena and analyticity for a two-component higher order
Camassa-Holm system,
Math. Nachr. (2018) 1--25. https://doi.org/10.1002/mana.201600469

\end{thebibliography}



\end{document}